\documentclass[11pt]{article}
\usepackage{amsmath, amssymb}
\usepackage{latexsym}
\usepackage{a4,epsfig}
\usepackage{amsfonts}
\usepackage{graphicx}
\usepackage{amscd}
\usepackage{amsbsy}

\numberwithin{equation}{section}
\newcommand\eq[1] {(\ref{#1})}

\newcommand{\bfm}[1]{\mbox{\boldmath ${#1}$}}

\newcommand{\beqa}{\begin{eqnarray}}
\newcommand{\eeqa}[1]{\label{#1}\end{eqnarray}}
\newcommand{\beq}{\begin{equation}}
\newcommand{\eeq}[1]{\label{#1}\end{equation}}
\newcommand{\Grad}{\nabla}

\newcommand{\Tr}{\mathop{\rm Tr}\nolimits}

\newcommand{\lang}{\langle}
\newcommand{\rang}{\rangle}
\newcommand{\Md}{\partial}

\newcommand{\Ga}{\alpha}
\newcommand{\Gb}{\beta}
\newcommand{\Gd}{\delta}
\newcommand{\Ge}{\epsilon}

\newcommand{\Gg}{\gamma}

\newcommand{\Gl}{\lambda}

\newcommand{\Gs}{\sigma}

\newcommand{\GD}{\Delta}

\newcommand{\GL}{\Lambda}

\newcommand{\GU}{\Upsilon}
\newcommand{\GO}{\Omega}

\newcommand{\BGs}{\bfm\sigma}

\def\Bc{{\bf c}}

\def\Be{{\bf e}}

\def\Bj{{\bf j}}
\def\Bk{{\bf k}}

\def\Bn{{\bf n}}

\def\Bq{{\bf q}}

\def\Bv{{\bf v}}

\def\Bx{{\bf x}}
\def\By{{\bf y}}

\def\BA{{\bf A}}
\def\BB{{\bf B}}

\def\BE{{\bf E}}

\def\BH{{\bf H}}
\def\BI{{\bf I}}
\def\BJ{{\bf J}}
\def\BK{{\bf K}}

\def\BP{{\bf P}}

\def\BR{{\bf R}}
\def\BS{{\bf S}}

\def\BU{{\bf U}}
\def\BV{{\bf V}}

\def \Tbb {{\mathbb T}}

\def \Lbb {{\mathbb L}}
\def \Ibb {{\mathbb I}}
\def \Abb {{\mathbb A}}
\def \Mbb {{\mathbb M}}

\def \dis {\displaystyle}

\def \p {\partial}
\def \Om {\Omega}

\def \ba {\begin{array}}
\def \ea {\end{array}}
\newtheorem {Thm} {Theorem} [section]

\newtheorem {Adef} [Thm] {Definition}

\newtheorem {Arem} [Thm] {Remark}

\newtheorem {Aexa} [Thm] {Example}

\newtheorem {Anot} [Thm] {Notation}

\def \refe #1.{(\ref{#1})}
\def \reff #1.{figure~\ref{#1}}
\def \refs #1.{section~\ref{#1}}
\def \refss #1.{subsection~\ref{#1}}
\def \refD #1.{Definition~\ref{#1}}
\def \refT #1.{Theorem~\ref{#1}}
\def \refL #1.{Lemma~\ref{#1}}
\def \refC #1.{Corollary~\ref{#1}}
\def \refP #1.{Proposition~\ref{#1}}
\def \refR #1.{Remark~\ref{#1}}
\def \refE #1.{Example~\ref{#1}}
\def \refN #1.{Notation~\ref{#1}}

\usepackage{graphics}
\usepackage{amssymb}
\begin{document}
\vspace{-1in}
\title{Bounds on the volume fractions of two materials in a three dimensional body from boundary measurements by the translation method}
\author{Hyeonbae Kang\\
\small{Department of Mathematics, Inha University, Incheon 402-751, Korea} \\
\small{(hbkang@inha.ac.kr)}\\
\\
Graeme W. Milton\\
\small{Department of Mathematics, University of Utah, Salt Lake City UT 84112, USA} \\
\small{(milton@math.utah.edu)}}
\date{}
\maketitle
\begin{abstract}
Using the translation method of Tartar, Murat, Lurie, and Cherkaev bounds are derived on the volume occupied by an inclusion
in a three-dimensional conducting body. They assume electrical impedance tomography measurements have been made 
for three sets of pairs of current flux and voltage measurements around the boundary. Additionally the conductivity of the 
inclusion and surrounding medium are assumed to be known. If the boundary data (Dirichlet or Neumann) is special, i.e. such that the fields inside the body 
would be uniform were the body homogeneous, then the bounds reduce to those of Milton and thus when the volume fraction is small to those
of Capdeboscq and Vogelius. 
\end{abstract}
\vskip2mm

\noindent Keywords: inverse problems, size estimation, electrical impedance tomography

\noindent
\vskip2mm
\section{Introduction}
\setcounter{equation}{0}
One of the most basic problems in electrical impedance tomography is to estimate the volume fraction of an
inclusion in a conducting body from one or more sets of electrical measurements of current flux and voltage pairs
around the boundary. If the boundary data (Dirichlet or Neumann) is special, i.e. such that the fields inside the body 
would be uniform were the body homogeneous, then sharp bounds on the volume fraction were derived by 
Nemat-Nasser and Hori \cite{Nemat-Nasser:1993:MOP, Nemat-Nasser:1995:UBO} (they did not present
their universal bounds as bounds on the volume fraction, but such bounds can be easily derived from
their results) and by Milton \cite{Milton:2011:UBE}. These reduce to the asymptotic bounds of 
Capdeboscq and Vogelius \cite{Capdeboscq:2003:OAE, Capdeboscq:2004:RSR} in the important case 
when the volume fraction is small. They are the natural generalization of the well-known 
Hashin-Shtrikman-Tartar-Murat-Lurie-Cherkaev bounds 
\cite{Hashin:1962:VAT, Hashin:1963:VAT, Murat:1985:CVH, Tartar:1985:EFC,  Lurie:1982:AEC, Lurie:1984:EEC}
for composites extended to two-phase bodies of arbitrary shape.

For other boundary conditions bounds, using the information from one measurement pair, 
were derived by Kang, Seo and Sheen \cite{Kang:1997:ICP}, Ikehata \cite{Ikehata:1998:SEI}, and 
Alessandrini, Morassi, Rosset, and Seo \cite{Alessandrini:1998:ICP, Alessandrini:2000:OSE, Alessandrini:2002:DCE}.
These other bounds involve constants which are not easy to determine. Another approach was taken by Kang, Kim and Milton
\cite{Kang:2011:SBV} who recognized that the translation method approach of Tartar, Murat, Lurie, and Cherkaev
\cite{Tartar:1979:ECH, Murat:1985:CVH, Tartar:1985:EFC, Lurie:1982:AEC, Lurie:1984:EEC, Milton:1990:CSP}
could be directly applied to bound the volume fraction using several measurement pairs, with arbitrary boundary conditions.
The constants in the resulting bounds can be easily evaluated from the boundary data. Furthermore the bounds 
were sharp if the boundary conditions and inclusion were such that the fields were uniform in the inclusion,
and shapes of inclusions having this property were found. For special boundary conditions they recovered the bounds 
of \cite{Milton:2011:UBE} and also provided new bounds, even
for this case. Those new bounds also reduced to the bounds of 
Capdeboscq and Vogelius \cite{Capdeboscq:2003:OAE, Capdeboscq:2004:RSR} for asymptotically small volume fractions. 
The same approach was applied to bound the size of an inclusion in a body for two-dimensional elasticity 
\cite{Milton:2011:BVF} (where another method, the method of splitting, was introduced)
and for the shallow shell equations \cite{Kang:2012:BVF}.

In this paper we extend the translation method approach to three-dimensional, two-phase, conducting bodies. This is relatively straightforward for
one of the bounds, as there is a null-Lagrangian associated with the electric fields, but requires additional care in the case of the
other bound as there is no null-Lagrangian associated with three-dimensional current fields and instead one uses quasiconvex
functions \cite{Tartar:1979:ECH}.

\section{Upper bounds}

\subsection{Derivation of the upper bounds}

Consider three potentials $V_1, V_2, V_3$ which satisfy $\Grad \cdot \Gs \Grad V_i=0$ in $\Om$, $i=1,2,3$, and $V_i=V_i^0$, and define $\Be_i =-\Grad V_i$ and $\Bj_i = \Gs \Be_i$. We rewrite this as
 \beq
 \Grad \cdot \Gs \Grad \BV=0 \quad\mbox{in }\Om,  \quad \BV=\BV^0 \quad \mbox{on } \p\Om,
 \eeq{UB1}
where $\BV=(V_1,V_2,V_3)$ is a vector potential, and following Tartar \cite{Tartar:1979:ECH} we define matrix valued fields
 \beq
 \BE =- \Grad \BV, \quad \BJ=\Gs \BE,
 \eeq{UB2}
where the notation is followed that $\BE$ has $\Be_1,\Be_2,\Be_3$ as columns (not rows) and $\BJ$ has $\Bj_1,\Bj_2,\Bj_3$ as columns.

We desire to use information about $q_1=-\Bj_1 \cdot \Bn$, $q_2=-\Bj_2 \cdot \Bn$, $q_3=-\Bj_3 \cdot \Bn$ and $V_1,V_2,V_3$ on $\p\Om$ to generate an upper bound on $f_1$ (the volume fraction of the phase with the highest conductivity).
The upper bound is expressed in terms of the response (or measurement) matrix $\BA$ defined by
 \beq
 \BA:= \lang \BE^T \Gs \BE \rang,
 \eeq{UB3}
which is computable from the boundary data:
 \beq
 \BA = \frac{1}{|\Om|} \int_\Om (-\Grad \BV)^T \BJ = \frac{1}{|\Om|} \int_{\p\Om} -\BV^T (\Bn \BJ) .
 \eeq{UB4}

We suppose that
 \beq
 \lang \BE \rang = \BI.
 \eeq{UB5}
We emphasize that this condition is achieved by taking a linear combination of the boundary data since
 \beq
 \lang \BE \rang =  \frac{1}{|\Om|} \int_{\p\Om} -\BV^T \Bn.
 \eeq{UB6}
We further assume that the response matrix $\BA$ is diagonal, namely,
 \beq
 \BA= \mbox{diag} \, [\Gl_1, \Gl_2, \Gl_3],
 \eeq{UB6-2}
by rotating the body $\GO$ and taking a further linear combination of $\BV^0$. In fact, since $\BA$ is symmetric, there is an orthogonal matrix $\BR$ such that $\BR^T \BA \BR$ is diagonal. We then take $\BR^T \BV^0 (\BR \Bx)$ as the new Dirichlet data. We emphasize that the condition \eq{UB5} still holds, since  $\lang \BE \rang$ transforms to $\BR^T\lang \BE \rang\BR=\BI$.

We introduce the 4th order tensor $\Tbb=(T_{ijkl})$ such that
 \beq
 \Tbb \BP = \Tr(\BP) \BI - \BP^T,
 \eeq{UB7}
where $\BP$ is a $3 \times 3$ matrix, {\it i.e.}, 
 \beq
 T_{ijkl} = \delta_{ij}\delta_{kl} - \delta_{il}\delta_{jk} .
 \eeq{UB7-1}
The tensor $\Tbb$ is associated with a null-Lagrangian and has the property that
$\Tbb (\Grad \BV)$ is divergence free, $\Grad \cdot \Tbb (\Grad \BV)=0$. In fact, since 
 \beq
 \Tbb (\Grad \BV)_{ij}=\delta_{ij} \frac{\p V_k}{\p x_k} - \frac{\p V_i}{\p x_j},
 \eeq{UB8}
we have
 \beq
 \frac{\p}{\p x_i} \Tbb (\Grad \BV)_{ij} = \frac{\p^2 V_k}{\p x_j \p x_k} - \frac{\p^2 V_i}{\p x_j \p x_i} =0.
 \eeq{UB9}
Note that the Einstein convention for summation is being used. If $\BV^{(1)}$ and $\BV^{(2)}$ are two solution vectors, then the matrix
 \beq
 \frac{1}{|\Om|} \int_{\Om} (\Grad \BV^{(1)})^T \Tbb (\Grad \BV^{(2)}) = \frac{1}{|\Om|} \int_{\p \Om} (\BV^{(1)})^T \Bn \Tbb (\Grad \BV^{(2)})
 \eeq{UB10}
can be computed from the boundary data as $\Bn  \Tbb (\Grad \BV^{(2)})$ only depends on the tangential derivatives of $\BV^{(2)}$ on $\p\Om$ ($\BV^{(2)}$ and $\Bn$ are row vectors). To see this it is enough to notice that the $i$-th component of the row vector $\Bn \Tbb (\Grad \BV)$ is
 \beq
 n_i (\Grad \cdot \BV) - \sum_j n_j \frac{\p V_j}{\p x_i} = \sum_{j} \left( n_i \frac{\p}{\p x_j} - n_j \frac{\p}{\p x_i} \right) V_j.
 \eeq{UB11}
We also note that
 \beq
 \frac{1}{|\Om|} \int_{\Om} (\Grad \BV^{(1)})^T \Tbb (\Grad \BV^{(2)}) = \frac{1}{|\Om|} \int_{\Om} (\Grad \BV^{(2)})^T \Tbb (\Grad \BV^{(1)}),
 \eeq{UB12}
since $\Tbb$ is symmetric in the sense that $T_{ijkl}=T_{klij}$ for all $i,j,k,l$.

Consider the quadratic form $\Tr \lang \BK^T \BE^T \Tbb (\BE  \BK)\rang$ on $3 \times 3$ matrices $\BK$. There is a unique symmetric 4th order tensor, say $\Mbb=(M_{ijkl})$, such that
 \beq
 \Tr \lang \BK^T \BE^T \Tbb (\BE  \BK)\rang = \BK: \Mbb \BK,
 \eeq{UB13}
where $\BA:\BB= \sum_{ij} A_{ij} B_{ij}$. $\Mbb$ being symmetric means $M_{ijkl}=M_{klij}$ for all $i,j,k,l$. Note that $M_{ijkl}$ is given by
 \beq
 M_{ijkl}= \Tr \lang (\BE\BI^{ij})^T \Tbb (\BE  \BI^{kl})\rang ,
 \eeq{UB14}
where $\BI^{ij}$ is the elementary matrix whose $(i,j)$ entry is 1 and other entries are $0$. We emphasize that $M_{ijkl}$ can be computed by the boundary data as was seen in \eq{UB10}. In fact, if we let $\BV^{ij}$ be the vector such that its $j$th component is $V_i$ and other components are $0$, then $(\Grad \BV)\BI^{ij} = \Grad \BV^{ij}$. Thus we have from \eq{UB10} and \eq{UB11} that 
 \begin{align*}
 \frac{1}{|\Om|} \int_{\Om} ((\Grad \BV) \BI^{ij})^T \Tbb ((\Grad \BV) \BI^{kl})
 &= \frac{1}{|\Om|} \int_{\Om} (\Grad \BV^{ij})^T \Tbb (\Grad \BV^{kl}) \\
 &= \frac{1}{|\Om|} \int_{\p \Om} \frac{\p V_k}{\p x_l} (\BV^{ij})^T \Bn - (\BV^{ij})^T \Bn (\Grad \BV^{kl})^T.
 \end{align*}
Therefore we get
 \beq
 M_{ijkl} = \frac{1}{|\Om|} \int_{\p \Om} V_i \left[ n_j \frac{\p V_k}{\p x_l} - n_l  \frac{\p V_k}{\p x_j} \right]. \eeq{UB16}
We may rewrite this in a symmetric form:
 \beq
 M_{ijkl} = \frac{1}{|\Om|} \int_{\Om} \left[ \frac{\p u_i}{\p x_j} \frac{\p u_k}{\p x_l} -\frac{\p u_i}{\p x_l} \frac{\p u_k}{\p x_j} \right],
 \eeq{UB17}
where $u_i$ is the solution to $\GD u_i=0$ in $\GO$ and $u_i=V_i^0$ on $\p\GO$, $i=1,2,3$. From this one can immediately see that $\Mbb$ is symmetric,
 \beq
 M_{ijkl} =0 \quad\mbox{if } i=k \ \mbox{or } j=l,
 \eeq{UB18}
and
 \beq
 M_{ijkl}= - M_{ilkj}.
 \eeq{UB19}
If affine boundary data are used, namely, $V_i^0=x_i$ for $i=1,2,3$, then one can easily see that $M_{1122}=M_{1133}=M_{2233}=1$ and hence that $\Mbb=\Tbb$.

For each $3 \times 3$ matrix $\BK$ and constant $c$ let
 \beq
 W_c := \Tr \lang \BK^T \BE^T (\Gs \Ibb + c \Tbb) (\BE  \BK)\rang
 = \Tr (\BK^T \BA \BK) + c \Tr \lang \BK^T \BE^T \Tbb (\BE  \BK)\rang,
 \eeq{UB20}
where $\Ibb$ is the 4th order identity tensor. Throughout this paper we identify $3 \times 3$ matrices with 9 dimensional vectors and 4th order tensors with $9 \times 9$ matrices following the usual convention:
 \begin{align*}
 & (1,1) \to 1, \ \ (2,1) \to 2, \ \ (3,1) \to 3, \ \ (1,2) \to 4, \ \ (2,2) \to 5, \\
 & (3,2) \to 6, \ \ (1,3) \to 7, \ \ (2,3) \to 8, \ \ (3,3) \to 9.
 \end{align*}
If we let $\Abb$ be the $9 \times 9$ matrix (4th order tensor) defined by
 \beq
 \Abb: = \begin{bmatrix} \BA & 0 & 0 \\
 0 & \BA & 0 \\
 0 & 0 & \BA
 \end{bmatrix},
 \eeq{UB21}
then $W_c$ can be written as
 \beq
 W_c =
 \begin{bmatrix} k_{11} \\ k_{21} \\ k_{31} \\ k_{12} \\ k_{22} \\ k_{32} \\ k_{13}  \\  k_{23} \\ k_{33} \end{bmatrix}
 \cdot (\Abb + c \Mbb)
 \begin{bmatrix} k_{11} \\ k_{21} \\ k_{31} \\ k_{12} \\ k_{22} \\ k_{32} \\ k_{13}  \\  k_{23} \\ k_{33} \end{bmatrix}.
 \eeq{UB22}

On the other hand, $W_c$ can be written as
 \beq
 W_c = \Tr \lang \BK^T \BE^T \Lbb_c (\BE \BK) \rang,
 \eeq{UB23}
where $\Lbb_c(\Bx)$ is the 4th order tensor $\Gs \Ibb + c \Tbb$ and the constant $c$ is chosen so that $\Lbb_c(\Bx)$ is positive definite for all $\Bx$. Since in the vector basis $\Tbb$ takes the form
 \beq
 \Tbb = \begin{bmatrix} 0 & 0 & 0 & 0 & 1 & 0 & 0 & 0 & 1 \\
 0 & 0 & 0 & -1 & 0 & 0 & 0 & 0 & 0 \\
 0 & 0 & 0 & 0 & 0 & 0 & -1 & 0 & 0 \\
 0 & -1 & 0 & 0 & 0 & 0 & 0 & 0 & 0 \\
 1 & 0 & 0 & 0 & 0 & 0 & 0 & 0 & 1 \\
 0 & 0 & 0 & 0 & 0 & 0 & 0 & -1 & 0 \\
 0 & 0 & -1 & 0 & 0 & 0 & 0 & 0 & 0 \\
 0 & 0 & 0 & 0 & 0 & -1 & 0 & 0 & 0 \\
 1 & 0 & 0 & 0 & 1 & 0 & 0 & 0 & 0 \end{bmatrix} ,
 \eeq{UB24}
we have (in the vector basis)
 \beq
 \Lbb_c= \begin{bmatrix} \Gs & 0 & 0 & 0 & c & 0 & 0 & 0 & c \\
 0 & \Gs & 0 & -c & 0 & 0 & 0 & 0 & 0 \\
 0 & 0 & \Gs & 0 & 0 & 0 & -c & 0 & 0 \\
 0 & -c & 0 & \Gs & 0 & 0 & 0 & 0 & 0 \\
 c & 0 & 0 & 0 & \Gs & 0 & 0 & 0 & c \\
 0 & 0 & 0 & 0 & 0 & \Gs & 0 & -c & 0 \\
 0 & 0 & -c & 0 & 0 & 0 & \Gs & 0 & 0 \\
 0 & 0 & 0 & 0 & 0 & -c & 0 & \Gs & 0 \\
 c & 0 & 0 & 0 & c & 0 & 0 & 0 & \Gs \end{bmatrix}.
 \eeq{UB25}

Since $\Tr \lang \BK^T \BE^T \Tbb (\BE  \BK)\rang$ is determined by the boundary data, the standard variational principle yields
 \begin{align}
 W_c &= \min_{\dis \underline{\BE} = -\Grad \underline{\BV} \atop \dis \underline{\BV} = \BV^0 \ \mbox{on } \p\Om}
 \lang \Tr (\BK^T \underline{\BE}^T \Gs (\underline{\BE} \BK)) \rang + c \Tr \lang \BK^T \BE^T \Tbb (\BE  \BK)\rang \nonumber \\
 &= \min_{\dis \underline{\BE} = -\Grad \underline{\BV} \atop \dis \underline{\BV} = \BV^0 \ \mbox{on } \p\Om}
 \lang \Tr (\BK^T \underline{\BE}^T \Lbb_c (\underline{\BE} \BK)) \rang. \label{UB26}
 \end{align}
We may rewrite this variational principle as
 \beq
 W_c= \min_{\dis \underline{\BE} = -\Grad \underline{\BV} \atop \dis \underline{\BV} = \BK^T\BV^0 \ \mbox{on } \p\Om}
 \lang \Tr (\underline{\BE}^T \Lbb_c \underline{\BE}) \rang.
 \eeq{UB27}
The constraints that $\lang \BE \rang = \BI$ and $\underline{\BV} = \BK^T\BV^0$ on $\p\Om$
imply  $\lang \underline{\BE} \rang = \lang \BE\BK \rang= \BK$, and so we have
 \beq
 W_c \ge \min_{\dis \underline{\BE} \atop \dis \lang \underline{\BE} \rang =\BK}
 \lang \Tr (\underline{\BE}^T \Lbb_c \underline{\BE}) \rang ,
 \eeq{UB28}
by relaxing the constraints.

One can easily see that the minimum of the right hand side of \eq{UB28} occurs when
 \beq
 \Lbb_c \underline{\BE} = \mu \ (\mbox{a constant matrix}),
 \eeq{UB29}
and the minimum is
 \beq
 \lang \Tr (\underline{\BE}^T \Lbb_c \underline{\BE}) \rang = \Tr (\BK^T \lang \Lbb_c^{-1} \rang^{-1} \BK ).
 \eeq{UB30}
Note that
 \beq
 \Tr (\BK^T \lang \Lbb_c^{-1} \rang^{-1} \BK ) = \sum_{i,j,k,l} k_{ij} (\lang \Lbb_c^{-1} \rang^{-1})_{ijkl} k_{kl}.
 \eeq{UB31}
So in the vector (rather than matrix) basis we get from \eq{UB22} and \eq{UB28}
 \beq
 \Abb + c \Mbb \ge \lang \Lbb_c^{-1} \rang^{-1}.
 \eeq{UB32}

By permuting rows and columns $\Lbb_c$ can be transformed to
 \beq
 \Lbb_c= \begin{bmatrix}
 \Gs & c & c &  &  &  &  &  &  \\
 c & \Gs & c &  &  &  &  &  &  \\
 c & c & \Gs &  &  &  &  &  &  \\
   &   &     & \Gs & -c &  &  &  &  \\
   &   &     & -c  & \Gs &  &  &  &  \\
   &   &     &     &     & \Gs & -c &  &   \\
   &   &     &     &    & -c  & \Gs &  &  \\
   &   &     &     &     &   &  & \Gs & -c \\
   &   &     &     &    &   &  & -c  & \Gs
   \end{bmatrix} .
 \eeq{UB33}
According to \eq{UB19}, $\Abb + c \Mbb$ takes (in that basis) the form
 \beq
 \begin{bmatrix} \Gl_1 & cM_{1122} & cM_{1133} &  &  &  &  &  &  \\
 cM_{1122} & \Gl_2 & cM_{2233} &  & * &  & * &  & * \\
 cM_{1133} & cM_{2233} & \Gl_3 &  &  &  &  &  &  \\
  &  * &  & \Gl_2 & -cM_{1122} &  & * &  & * \\
  &  &  & -cM_{1122} & \Gl_1 &  &  &  &  \\
  &  *&  &  & * & \Gl_3 & -cM_{1133} &  & * \\
  &  & & &  & -cM_{1133} & \Gl_1 &  & \\
  & *& &  & * &  & * & \Gl_3 & -cM_{2233} \\
  &  &  &  &  &  &  & -cM_{2233} & \Gl_2 \end{bmatrix}.
 \eeq{UB34}
It is worth emphasizing that $\Lbb_c^{-1}$ takes the same block form as $\Lbb_c$. By comparing $3 \times 3$ blocks on the left upper corner, \eq{UB32} yields
 \beq
 \begin{bmatrix} \Gl_1 & cM_{1122} & cM_{1133} \\
 cM_{1122} & \Gl_2 & cM_{2233} \\
 cM_{1133} & cM_{2233} & \Gl_3
 \end{bmatrix} \ge \left\lang \begin{bmatrix}
 \Gs & c & c \\
 c & \Gs & c \\
 c & c & \Gs
 \end{bmatrix}^{-1} \right\rangle^{-1} .
 \eeq{UB35}
Note that the inverse of $\begin{bmatrix} \Gs & c & c \\ c & \Gs & c \\ c & c & \Gs \end{bmatrix}$ is of the form $\begin{bmatrix} a & b & b \\ b & a & b \\ b & b & a \end{bmatrix}$ where
 \beq
 a= - \frac{b(\Gs+c)}{c} \quad\mbox{and}\quad b= \frac{c}{(c-\Gs)(2c+\Gs)},
 \eeq{UB36}
satisfy
 \beq
 a+2b=\frac{1}{\Gs+2c}.
 \eeq{UB37}
Thus we have
 \beq
 \left\langle \begin{bmatrix}
 \Gs & c & c \\
 c & \Gs & c \\
 c & c & \Gs
 \end{bmatrix}^{-1} \right\rangle^{-1} = \begin{bmatrix}
 \Ga & \Gb & \Gb \\
 \Gb & \Ga & \Gb \\
 \Gb & \Gb & \Ga
 \end{bmatrix},
 \eeq{UB38}
where
 \beq
 \Ga= - \frac{\Gb (\langle a \rangle + \langle b \rangle)}{\langle b \rangle} \quad\mbox{and}\quad
 \Gb= \frac{\langle b \rangle}{(\langle b \rangle - \langle a \rangle) (2\langle b \rangle + \langle a \rangle)}.
 \eeq{UB39}
One can then see that as $c$ is increasing and approaching $\Gs_2$,  $ \left\langle \begin{bmatrix}
 \Gs & c & c \\
 c & \Gs & c \\
 c & c & \Gs
 \end{bmatrix}^{-1} \right\rangle^{-1}$ is decreasing (in terms of eigenvalues) and approaching 
$\begin{bmatrix}
 d & d & d  \\
 d & d & d  \\
 d & d & d  \end{bmatrix}$,
where
 \beq
 d= \left(\frac{f_2}{\Gs_2}+ \frac{3f_1}{2\Gs_2+\Gs_1} \right)^{-1}>\Gs_2 .
 \eeq{UB40}
This shows in particular that $ \left\langle \begin{bmatrix}
 \Gs & c & c \\
 c & \Gs & c \\
 c & c & \Gs
 \end{bmatrix}^{-1} \right\rangle^{-1} - \Gs_2 \begin{bmatrix}
 1 & 1 & 1  \\
 1 & 1 & 1  \\
 1 & 1 & 1  \end{bmatrix}$ is positive-definite if $c<\Gs_2$ regardless of the volume fraction $f_1$. It then follows from \eq{UB35} that
 \beq
 \begin{bmatrix} \Gl_1 & cM_{1122} & cM_{1133} \\
 cM_{1122} & \Gl_2 & cM_{2233} \\
 cM_{1133} & cM_{2233} & \Gl_3
 \end{bmatrix} - \Gs_2 \begin{bmatrix}
 1 & 1 & 1  \\
 1 & 1 & 1  \\
 1 & 1 & 1  \end{bmatrix} \ge \left\lang \begin{bmatrix}
 \Gs & c & c \\
 c & \Gs & c \\
 c & c & \Gs
 \end{bmatrix}^{-1} \right\rangle^{-1} - \Gs_2 \begin{bmatrix}
 1 & 1 & 1  \\
 1 & 1 & 1  \\
 1 & 1 & 1  \end{bmatrix} ,
 \eeq{UB41}
or
 \beq
 \Bv\cdot\left(\begin{bmatrix} \Gl_1 & cM_{1122} & cM_{1133} \\
 cM_{1122} & \Gl_2 & cM_{2233} \\
 cM_{1133} & cM_{2233} & \Gl_3
 \end{bmatrix} - \Gs_2 \begin{bmatrix}
 1 & 1 & 1  \\
 1 & 1 & 1  \\
 1 & 1 & 1  \end{bmatrix} \right)^{-1}\Bv \le \Bv\cdot\left( \left\lang \begin{bmatrix}
 \Gs & c & c \\
 c & \Gs & c \\
 c & c & \Gs
 \end{bmatrix}^{-1} \right\rangle^{-1} - \Gs_2 \begin{bmatrix}
 1 & 1 & 1  \\
 1 & 1 & 1  \\
 1 & 1 & 1  \end{bmatrix} \right)^{-1}\Bv .
 \eeq{UB42}
Note that
 \beq
 \left( \left\lang \begin{bmatrix}
 \Gs & c & c \\
 c & \Gs & c \\
 c & c & \Gs
 \end{bmatrix}^{-1} \right\rangle^{-1} - \Gs_2 \begin{bmatrix}
 1 & 1 & 1  \\
 1 & 1 & 1  \\
 1 & 1 & 1  \end{bmatrix} \right)^{-1} =
 \left( \begin{bmatrix}
 \Ga & \Gb & \Gb \\
 \Gb & \Ga & \Gb \\
 \Gb & \Gb & \Ga
 \end{bmatrix} - \Gs_2 \begin{bmatrix}
 1 & 1 & 1  \\
 1 & 1 & 1  \\
 1 & 1 & 1  \end{bmatrix} \right)^{-1} = \begin{bmatrix}
 \Gg & \Gd & \Gd \\
 \Gd & \Gg & \Gd \\
 \Gd & \Gd & \Gg
 \end{bmatrix},
 \eeq{UB43}
where
 \beq
 \Gg= - \frac{\Gd (\Ga-\Gs_2 + \Gb-\Gs_2)}{\Gb-\Gs_2}, \quad \Gd= \frac{\Gb-\Gs_2}{(\Gb-\Ga)(2\Gb-2\Gs_2+\Gs-\Gs_2)}.
 \eeq{UB44}
Note that $\Gd$ blows up as $c$ approaches $\Gs_2$, in phase 2. Thus to get a non-trivial bound we take
$\Bv= [1 \ 1 \ 1]^T$ and we have
 \beq
 \begin{bmatrix} 1 \\ 1 \\ 1 \end{bmatrix} \cdot \begin{bmatrix} \Gg & \Gd & \Gd \\  \Gd & \Gg & \Gd \\
 \Gd & \Gd & \Gg \end{bmatrix}
 \begin{bmatrix} 1 \\ 1 \\ 1 \end{bmatrix} = 3(\Gg+2\Gd)= \frac{3}{\frac{1}{\left\langle \frac{1}{\Gs+2c} \right\rangle} - 3\Gs_2} \to \frac{1}{d-\Gs_2},
 \eeq{UB45}
as $c \to \Gs_2$ where $d$ is given by \eq{UB40}.

The left hand side of \eq{UB42} approaches $T$ as $c \to \Gs_2$ where
 \beq
 T:= \begin{bmatrix} 1 \\ 1 \\ 1 \end{bmatrix} \cdot
 \left(\begin{bmatrix} \Gl_1-\Gs_2 & 0 & 0 \\
 0 & \Gl_2-\Gs_2 & 0 \\
 0 & 0 & \Gl_3-\Gs_2
 \end{bmatrix} - \Gs_2 \begin{bmatrix} 0 & M_{1122}-1 & M_{1133}-1 \\
 M_{1122}-1 & 0 & M_{2233}-1 \\
 M_{1133}-1 & M_{2233}-1 & 0
 \end{bmatrix} \right)^{-1} \begin{bmatrix} 1 \\ 1 \\ 1 \end{bmatrix} .
 \eeq{UB46}
Thus we have the following upper bound:
 \beq
 f_1 \le \frac{\Gs_1 + 2 \Gs_2}{\Gs_1-\Gs_2} \frac{1}{1+ \Gs_2 T}.
 \eeq{UB47}
This is generally not the tightest bound that can be obtained from \eq{UB32} as we have ignored the off-diagonal elements marked * in 
\eq{UB34}. A tighter bound (but one which is not so easily solved as an inequality on $f_1$) is provided by the appropriate interval 
of volume fractions $f_1$ given by the inequality (implied by \eq{UB32})
\beq \det[\Abb + \Gs_2 \Mbb- \lang \Lbb_{\Gs_2}^{-1} \rang^{-1}]\geq 0,
\eeq{UB47.a}
in which, because $\Lbb_{\Gs_2}$ is singular, $\lang \Lbb_{\Gs_2}^{-1} \rang^{-1}$ is defined as $\lim_{c \to \Gs_2} \lang \Lbb_{c}^{-1} \rang^{-1}$.
In fact, one can show by tedious computations that in the basis in which $\Lbb_c$ is given by \eq{UB33}
 \beq
 \lang \Lbb_{\Gs_2}^{-1} \rang^{-1} = \begin{bmatrix}
 a & a & a &  &  &  &  &  &  \\
 a & a & a &  &  &  &  &  &  \\
 a & a & a &  &  &  &  &  &  \\
   &   &     & b & -b &  &  &  &  \\
   &   &     & -b  & b &  &  &  &  \\
   &   &     &     &     & b & -b &  &   \\
   &   &     &     &    & -b  & b &  &  \\
   &   &     &     &     &   &  & b & -b \\
   &   &     &     &    &   &  & -b  & b
   \end{bmatrix} ,
 \eeq{AB00}
where
 \beq
 a= \left(\frac{f_2}{\Gs_2}+ \frac{3f_1}{2\Gs_2+\Gs_1} \right)^{-1}, \quad b=\frac{\Gs_2}{f_1}.
 \eeq{AB00-1}

Note also that the translation tensor $\Tbb$ is not the most general one: we could have taken $\Tbb$ to be any
symmetric fourth-order tensor satisfying the symmetries $T_{ijkl}= - T_{kjil}$. We have not explored this freedom in the choice $\Tbb$,
which could lead to even tighter bounds.

\subsection{Special boundary conditions}

If the affine Dirichlet boundary conditions are prescribed, namely, $\BV^0=\Bx$ on $\p\Om$, then $\lang \BE \rang = \BI$, and hence
 \beq
 \BA= \lang \BE^T \Gs \BE \rang = \lang \BE^T \rang \lang \Gs \BE \rang = \BI \BGs_D \lang \BE \rang =\BGs_D,
 \eeq{UB48}
where $\BGs_D$ is the Dirichlet tensor. Since $M_{1122}=M_{1133}=M_{2233}=1$ in this case, we have
 \beq
 T= \Tr [(\BGs_D - \Gs_2 \BI)^{-1}],
 \eeq{UB49}
and hence
 \beq
 f_1 \le \frac{\Gs_1 + 2 \Gs_2}{\Gs_1-\Gs_2} \frac{1}{1+ \Gs_2 \Tr [(\BGs_D - \Gs_2 \BI)^{-1}]}.
 \eeq{UB50}
This upper bound is exactly the same as the one obtained in \cite{Milton:2011:UBE}.

The inequality \eq{UB32} yields not only \eq{UB35} but also
 \beq
 \begin{bmatrix}
 \Gl_2 & -cM_{1122} \\
 -cM_{1122} & \Gl_1 \end{bmatrix}
 \ge \left\lang \begin{bmatrix}
 \Gs & -c \\
 -c  & \Gs
 \end{bmatrix}^{-1} \right\rangle^{-1} ,
 \eeq{UB51}
 \beq
 \begin{bmatrix}
 \Gl_3 & -cM_{1133} \\
 -cM_{1133} & \Gl_1 \end{bmatrix}
 \ge \left\lang \begin{bmatrix}
 \Gs & -c \\
 -c  & \Gs
 \end{bmatrix}^{-1} \right\rangle^{-1} ,
 \eeq{UB52}
and
 \beq
 \begin{bmatrix} \Gl_3 & -cM_{2233} \\
 -cM_{2233} & \Gl_2 \end{bmatrix}
 \ge \left\lang \begin{bmatrix}
 \Gs & -c \\
 -c  & \Gs
 \end{bmatrix}^{-1} \right\rangle^{-1}.
 \eeq{UB53}
When $\BV^0=\Bx$ on $\p\Om$ these inequalities do not yield a better bound than \eq{UB47}. To see this, observe that with these boundary conditions \eq{UB51} becomes
 \beq
 \begin{bmatrix}
 \Gl_2 & -c \\
 -c & \Gl_1 \end{bmatrix}
 \ge \left\lang \begin{bmatrix}
 \Gs & -c \\
 -c  & \Gs
 \end{bmatrix}^{-1} \right\rangle^{-1} ,
 \eeq{UB54}
and hence
 \beq
 \left( \begin{bmatrix}
 \Gl_2 & -c \\
 -c & \Gl_1 \end{bmatrix} - \Gs_2 \begin{bmatrix}
 1 & -1 \\
 -1 & 1 \end{bmatrix} \right)^{-1}
 \ge \left( \left\lang \begin{bmatrix}
 \Gs & -c \\
 -c  & \Gs
 \end{bmatrix}^{-1} \right\rangle^{-1} - \begin{bmatrix}
 1 & -1 \\
 -1 & 1 \end{bmatrix} \right)^{-1}.
 \eeq{UB55}
We then obtain
 \beq
 \frac{1}{\Gl_1-\Gs_2} + \frac{1}{\Gl_2-\Gs_2} \le \frac{1}{\frac{1}{e} - \Gs_2},
 \eeq{UB56}
where
 \beq
 e=\left(\frac{f_2}{2\Gs_2}+ \frac{f_1}{\Gs_2+\Gs_1} \right)^{-1} .
 \eeq{UB57}
So far we have shown that \eq{UB51}, \eq{UB52} and \eq{UB53} yield
 \beq
 \frac{1}{\Gl_1-\Gs_2} + \frac{1}{\Gl_2-\Gs_2} + \frac{1}{\Gl_3-\Gs_2} \le \frac{\frac{3}{2}}{\frac{1}{e} - \Gs_2},
 \eeq{UB58}
which in turn yields
 \beq
 f_1 \le \frac{\Gs_1 + \Gs_2}{\Gs_1-\Gs_2} \frac{1}{1+ \frac{2}{3}\Gs_2 \Tr [(\BGs_D - \Gs_2 \BI)^{-1}]}.
 \eeq{UB59}

Note that
 \beq
 \Tr [(\BGs_D - \Gs_2 \BI)^{-1}] \ge \frac{3}{\Gs_1-\Gs_2}.
 \eeq{UB60}
Using this fact, one can easily show that
 \beq
 \frac{\Gs_1 + 2 \Gs_2}{1+ \Gs_2 \Tr [(\BGs_D - \Gs_2 \BI)^{-1}]} \le
 \frac{\Gs_1 + \Gs_2}{1+ \frac{2}{3}\Gs_2 \Tr [(\BGs_D - \Gs_2 \BI)^{-1}]} .
 \eeq{UB61}
Thus the bound in \eq{UB50} is better than the one in \eq{UB59}. They are the same bounds only when $f_1=1$.

\section{Attainability of the determinant bound \eq{UB47.a}}
Suppose that the bound in \eq{UB28} is attained for a field $\BE$. Then one can show (see \cite{Kang:2011:SBV} for more details in the analogous two-dimensional case)
that
 \beq
 \Lbb_{\Gs_2} \BE = \lang \Lbb_{\Gs_2}^{-1} \rang^{-1} \lang \BE \rang,
 \eeq{AB0}
where $\lang \Lbb_{\Gs_2}^{-1} \rang^{-1}$, given by \eq{AB00}, is defined as $\lim_{c \to \Gs_2} \lang \Lbb_{c}^{-1} \rang^{-1}$.
We now show that the attainability condition \eq{AB0} holds if and only if
 \beq
 \Lbb_{\Gs_2} \BE = \BU,
 \eeq{AB00-2}
for some constant matrix $\BU$. The only if part is trivial. The converse is clear heuristically:
Since $\BE = \Lbb_{\Gs_2}^{-1} \BU$ and $\BU$ is constant, $\lang \BE \rang = \lang \Lbb_{\Gs_2}^{-1} \rang \BU$ and hence $\BU = \lang \Lbb_{\Gs_2}^{-1} \rang^{-1} \lang \BE \rang$, therefore \eq{AB0} holds. To show the converse rigorously, we use the basis in which $\Lbb_c$ is given by \eq{UB33}. In that basis,
 \beq
 \BE=\begin{bmatrix} E_{11} & E_{22} & E_{33} & E_{21} & E_{12} & E_{31} & E_{13} & E_{32} & E_{23} \end{bmatrix}^T.
 \eeq{AB00-3}
Let $\Lbb_{\Gs_2}^{(1)}$ and $\Lbb_{\Gs_2}^{(2)}$ denote $\Lbb_{\Gs_2}$ in phase 1 and 2, respectively, and likewise $\BE^{(1)}$ and $\BE^{(2)}$. Since $\BU$ is in the image of $\Lbb_{\Gs_2}^{(2)}$, $\BU$ takes the form
 \beq
 \BU=\begin{bmatrix} \Ga & \Ga & \Ga & \Gb_1 & -\Gb_1 & \Gb_2 & -\Gb_2 & \Gb_3 & -\Gb_3 \end{bmatrix}^T.
 \eeq{AB00-4}
Thus $\BE^{(2)}$ takes the form
 \beq
 \BE^{(2)}=\BE_0 + \BE_{\perp},
 \eeq{AB00-5}
where
 \beq
 \BE_0=\begin{bmatrix} \frac{\Ga}{3\Gs_2} & \frac{\Ga}{3\Gs_2} & \frac{\Ga}{3\Gs_2} & \frac{\Gb_1}{2\Gs_2} & -\frac{\Gb_1}{2\Gs_2} & \frac{\Gb_2}{2\Gs_2} & -\frac{\Gb_2}{2\Gs_2} & \frac{\Gb_3}{2\Gs_2} & -\frac{\Gb_3}{2\Gs_2} \end{bmatrix}^T,
 \eeq{AB00-6}
and $\BE_{\perp}$ is in the kernel of $\Lbb_{\Gs_2}^{(2)}$, {\it i.e.}, $\BE_{\perp}$ is symmetric and $\Tr \BE_{\perp}=0$. Since $\int_{\mbox{\small phase 2}} \BE_{\perp}$ is also symmetric and of trace $0$, we can see from the form \eq{AB00} of $\lang \Lbb_{\Gs_2}^{-1} \rang^{-1}$ that
 \beq
 \lang \Lbb_{\Gs_2}^{-1} \rang^{-1} \int_{\mbox{\small phase 2}} \BE_{\perp} =0.
 \eeq{AB00-7}
Thus we obtain
 \beq
 \lang \Lbb_{\Gs_2}^{-1} \rang^{-1} \lang \BE \rang = \lang \Lbb_{\Gs_2}^{-1} \rang^{-1} ( f_1 (\Lbb_{\Gs_2}^{(1)})^{-1} \BU + f_2 \BE_0 ) = \BU,
 \eeq{AB00-8}
where the last equality can be shown by an elementary but tedious computation.

The condition \eq{AB00-2} implies in particular that $\BE$ is constant in phase 1. We now show that this condition alone guarantees that the upper bound is attained.

Suppose that $\BE$ is constant in phase 1. After coordinate changes and taking linear combinations of potentials if necessary, we may assume that
 \beq
 \BE= \Grad \BV=\BI \quad\mbox{in phase 1}.
 \eeq{AB01}
We claim that if \eq{AB01} holds, then there exists a potential $\psi$ in phase 2 such that
 \beq
 \BV = \Grad \psi \quad \mbox{and} \quad \Delta \psi = \frac{\Gs_1+2\Gs_2}{\Gs_2} \quad \mbox{in phase 2}.
 \eeq{AB1}

To show \eq{AB1} let $D$ be a connected component of phase 1 and $\Om_0$ an open subset of $\Om$ containing $D$ such that $\Om_0 \setminus D$ is phase 2 and connected. 
It is known (see section 2.5 of \cite{Ammari:2004:RSI}) that for $j=1,2,3$, there is a harmonic function
$H_j$ and a potential $\varphi_j$ such that
 \beq
 V_j (\Bx) = H_j(\Bx) + \int_{\p D} \frac{\varphi_j(\By)}{4\pi |\Bx- \By|} d\sigma(\By), \quad \Bx
 \in \Om_0.
 \eeq{AB2}
Furthermore, we have
 \beq
 \varphi_j = \frac{\Gs_1-\Gs_2}{\Gs_2} \frac{\p V_j}{\p\Bn} \Big|_{-} \quad\mbox{on } \p D.
 \eeq{AB3}
If \eq{AB01} holds, then
 \beq
 \varphi_j = \frac{\Gs_1-\Gs_2}{\Gs_2} n_j  \quad\mbox{on } \p D,
 \eeq{AB4}
and hence
 \beq
 V_j (\Bx) = H_j(\Bx) - \frac{\Gs_1-\Gs_2}{\Gs_2} \frac{\p}{\p x_j} \int_{D} \frac{1}{4\pi |\Bx- \By|} d\By, \quad \Bx \in \Om_0.
 \eeq{AB5}
Thus we have
 \beq
 H_j(\Bx) = x_j + c_j + \frac{\Gs_1-\Gs_2}{\Gs_2} \frac{\p}{\p x_j} \int_{D} \frac{1}{4\pi |\Bx- \By|} d\By, \quad \Bx \in D,
 \eeq{AB6}
for a constant $c_j$. Let $\BH=(H_1, H_2, H_3)$ and $\Bc=(c_1, c_2, c_3)$. Then,
 \beq
 \BH(\Bx) = \Grad \left( \frac{1}{2} |\Bx|^2 + \Bc \cdot \Bx + \frac{\Gs_1-\Gs_2}{\Gs_2} \int_{D} \frac{1}{4\pi |\Bx- \By|} d\By \right), \quad \Bx \in D.
 \eeq{AB7}
In particular, we have $\Grad \times \BH=0$ in $D$. Further we have $\Grad \cdot \BH= \frac{\Gs_1+2\Gs_2}{\Gs_2}$ in $D$. Since all components of $\Grad \times \BH$ and $\Grad \cdot \BH$ are harmonic in $\Om_0$, we have $\Grad \times \BH=0$ and $\Grad \cdot \BH= \frac{\Gs_1}{\Gs_2}$ in $\Om_0$. It then follows from \eq{AB5} that $\Grad \times \BV=0$ and $\Grad \cdot \BV= \frac{\Gs_1}{\Gs_2}$ in $\Om_0 \setminus \overline{D}$, and hence in phase 2 since the phase 2 is connected.

This does not quite guarantee existence of a potential $\psi$ such that $\BV = \Grad \psi$ since phase 2 is not simply connected. But, $\BV$ can be written as $\BV=\BH+\BS$ where $\BH$ is harmonic in $\Om$ and $\BS$ is harmonic in $\mathbb{R}^3 \setminus \{\mbox{phase 1} \}$ and $\BS(\Bx) \to 0$ as $|\Bx| \to \infty$. What we have shown above is that $\Grad \times \BH=0$ in $\Om$ and $\Grad \times \BS=0$ in $\mathbb{R}^3 \setminus \{\mbox{phase 1} \}$. It implies that there are potentials $\psi_1$ and $\psi_2$ such that $\BH= \Grad \psi_1$ in $\Om$ and $\BS= \Grad \psi_2$ in $\mathbb{R}^3 \setminus \{\mbox{phase 1} \}$. Thus we have \eq{AB1}.

The relation \eq{AB1} implies that $\BE$ is symmetric and $\Tr \BE$ is constant in phase 2. Therefore $\Lbb_{\Gs_2} \BE$ is constant in phase 2. Because of the continuity of $\Gs \BE$ along the interface, 
it follows that $\Lbb_{\Gs_2} \BE$ is constant in $\Om$ and hence the determinant bound \eq{UB47.a} is attained.

\section{Lower bounds}

In this section we derive lower bounds for $f_1$. We still use the translation method. But as there is no null-Lagrangian for current fields (see below) 
we use a quasiconvex quadratic form, namely the extremal one discussed on page 546 of \cite{Milton:2002:TOC}.

Let $\GL_h$, $\GL_s$ and $\GL_a$ be orthogonal projections onto the space of $3 \times 3$ matrices proportional to $\BI$, symmetric and trace free, and antisymmetric, respectively, namely
 \begin{align*}
 \GL_h \BP &= \frac{1}{3} \Tr(\BP) \BI, \\
 \GL_s \BP &= \frac{1}{2} (\BP+\BP^T) - \GL_h \BP, \\
 \GL_a \BP &= \frac{1}{2} (\BP-\BP^T),
 \end{align*}
for any $3 \times 3$ matrix $\BP$. Define a 4th order tensor $\Tbb^\prime$ by
 \beq
 \Tbb^\prime : = 2\GL_s - \GL_h,~~~{\rm so~that~}\Tbb^\prime\BP=\BP+\BP^T- \Tr(\BP) \BI.
 \eeq{LB1}
This tensor is positive semidefinite on the set of $3 \times 3$ rank two matrix valued fields, and isotropic in the sense that
 \beq
 \Tbb^\prime (\BR^T \BP \BR) =\BR^T (\Tbb^\prime \BP)\BR,
 \eeq{LB1-4}
for any unitary transform $\BR$. Due to this isotropy it suffices, to check for positive semidefiniteness on rank two matrices, to consider matrices of the form 
\beq \BP=\begin{bmatrix}
p_{11} & p_{12} & p_{13} \\
p_{21} & p_{22} & p_{23} \\
0 & 0 & 0 
\end{bmatrix},
\eeq{LB1b}
and in that case
\beq  \Tr (\BP^T \Tbb^\prime \BP)=(p_{11}-p_{22})^2+(p_{12}+p_{21})^2+p_{13}^2+p_{23}^2
\eeq{LB1c}
is clearly non-negative and zero only when $p_{11}=p_{22}$, $p_{12}=-p_{21}$, and $p_{13}=p_{23}=0$. More generally, $\Tr (\BP^T \Tbb^\prime \BP)$ will be zero if and only if
the rank two matrix $\BP$ is of the form
\beq  P_{ij}=\Ga_0(k_ik_j-\Gd_{ij}\Bk\cdot\Bk)+\Gb_0\Ge_{ijm}k_m, \eeq{LB1ca}
for some constants $\Ga_0$ and $\Gb_0$, and for some vector $\Bk$ which is a vector such that $\BP\Bk=0$. Here $\Ge_{klm}$ is the completely antisymmetric 
Levi-Civita symbol, taking value $1$ if $klm$ is an even permutation of $123$, $-1$ if it is an odd permutation, and zero otherwise.

If $\underline{\BJ}(\Bx)$ is periodic, with unit cell $\GU$, 
and has zero divergence then its Fourier components
$\underline{\hat\BJ}(\Bk)$ satisfy $\Bk\cdot\underline{\hat\BJ}(\Bk)=0$ and hence for $\Bk\ne 0$ are matrices of rank two or less. As observed by Murat and Tartar 
\cite{Tartar:1979:CCA, Murat:1985:CVH, Tartar:1985:EFC}, it then 
follows by a straightforward application of Plancherel's theorem that the quadratic form associated with $\Tbb^\prime$ is quasiconvex on divergence free periodic fields:
\beq \lang\Tr (\underline{\BJ}^T \Tbb^\prime\underline{\BJ})\rang_{\GU}\geq \Tr (\lang\underline{\BJ}^T\rang_{\GU} \Tbb^\prime\lang\underline{\BJ}\rang_{\GU}),
\eeq{LB1d}
(where the angular brackets $\lang\cdot\rang_{\GU}$ denote volume averages over $\GU$) with equality if and only if the Fourier components of $\underline{\BJ}$ have
the form \eq{LB1ca}, {\it i.e.}
\beq \underline{J}_{kl}=J^0_{kl}+\frac{\Md^2\Ga}{\Md x_k \Md x_l}-\Gd_{kl}\Delta\Ga+\Ge_{klm}\frac{\Md\Gb}{\Md x_m}
\eeq{LB1e}
for some scalar potentials $\Ga(\Bx)$ and $\Gb(\Bx)$, where the $J^0_{kl}$ are the elements of a constant matrix $\BJ^0$. Note that $\underline{\BJ}$ given by \eq{LB1e} is divergence free, has $\BJ^0$ as its average over $\GU$,
and 
\beq  \Tbb^\prime\underline{\BJ}=\Tbb^\prime\BJ^0+2\nabla\nabla\Ga \eeq{LB1f}
is the gradient of a potential.

Incidentally, there is no symmetric fourth order tensor $\Tbb^\prime$, with elements $T_{ijk\ell}^\prime$, such that $\Tbb^\prime$ 
is a null Lagrangian for currents. If there were the quadratic form $\Tr (\BP^T \Tbb^\prime \BP)$ would be zero for all $3 \times 3$ rank two matrix valued fields. In
particular, being zero for all matrices of the form \eq{LB1b} implies $T_{ijk\ell}^\prime=0$ unless $i=3$ or $k=3$. Similarly, by considering matrices where the
first or second row of $\BP$ is zero we deduce that $T_{ijk\ell}^\prime=0$ unless $i=1$ or $k=1$, and $T_{ijk\ell}^\prime=0$ unless $i=2$ or $k=2$. It then follows
that all elements of $\Tbb^\prime$ must be zero. 

Let $\Bj_1$, $\Bj_2$ and $\Bj_3$ be the currents corresponding to the measurements $\Bj_1 \cdot \Bn=q_1$, $\Bj_2 \cdot \Bn=q_2$ and $\Bj_3 \cdot \Bn=q_3$ on $\p\GO$. Let $\BJ = [ \Bj_1 \ \Bj_2 \ \Bj_3 ]$ and suppose that
 \beq
 \lang \BJ \rang = \BI.
 \eeq{LB2}
Define the response matrix $\BA^\prime$ by
 \beq
 \BA^\prime:= \lang \BJ^T \Gs^{-1} \BJ \rang,
 \eeq{LB3}
which is computable from the boundary data:
 \beq
 \BA^\prime = \frac{1}{|\Om|} \int_\Om \BJ^T (-\Grad \BV) = \frac{1}{|\Om|} \int_{\p\Om} -(\Bn \BJ)^T \BV .
 \eeq{LB4}

For a positive constant $c$, define
 \beq
 \Lbb_c^\prime: = \Gs^{-1} \Ibb - c \Tbb^\prime= (\Gs^{-1}+c) \GL_h + (\Gs^{-1}-2c) \GL_s + \Gs^{-1}\GL_a.
 \eeq{LB5}
We assume that
 \beq
 c < \Gs_1^{-1}/2,
 \eeq{LB6}
so that $\Lbb_c^\prime$ is positive definite. It is helpful to recall that $\Gs_2 < \Gs_1$.

Define
 \beq
 W_c^\prime := \min_{\nabla \cdot \underline{\BJ} = 0, \ \lang \underline{\BJ} \rang = \BI, \atop \Bn \underline{\BJ} = -\Bq \ \textrm{on } \p\Om}
 \Tr \lang \underline{\BJ}^T \Lbb_c^\prime \underline{\BJ} \rang,
 \eeq{LB10}
where $\Bq=[q_1 \ q_2 \ q_3]$. Let
 \beq
 g:= \min_{\nabla \cdot \underline{\BJ} = 0, \ \lang \underline{\BJ} \rang = \BI, \atop \Bn \underline{\BJ} = -\Bq \ \textrm{on } \p\Om}
 \Tr \lang \underline{\BJ}^T \Tbb^\prime \underline{\BJ} \rang .
 \eeq{LB10-2}
Then we have
 \beq
 W_c^\prime \le \min_{\nabla \cdot \underline{\BJ} = 0, \ \lang \underline{\BJ} \rang = \BI, \atop \Bn \underline{\BJ} = -\Bq \ \textrm{on } \p\Om}
 \Tr \lang \underline{\BJ}^T \Gs^{-1} \underline{\BJ} \rang - \min_{\nabla \cdot \underline{\BJ} = 0, \ \lang \underline{\BJ} \rang = \BI, \atop \Bn \underline{\BJ} = -\Bq \ \textrm{on } \p\Om}
 c\Tr \lang \underline{\BJ}^T \Tbb^\prime \underline{\BJ} \rang = \Tr \BA^\prime - cg \le \Tr \BA^\prime - cg^-,
 \eeq{LB10-4}
where $g^-$ is a lower bound on $g$. Lower bounds on $g$ are easily obtained. Choose a cube $\GU$ which contains $\GO$ and construct a divergence free field $\underline{\BJ}$ in $\GU\setminus\GO$
such that $\underline{\BJ}$ satisfies periodic boundary conditions on $\p\GU$ and $\Bn \underline{\BJ} = -\Bq$ on $\p\Om$. Then from \eq{LB1d} the inequality $g^-\le g$ holds with
\beq g^-=\frac{1}{p}\left\{\Tr[(p\BI+(1-p)\lang\underline{\BJ}\rang_{\GU\setminus\GO})^T\Tbb^\prime(p\BI+(1-p)\lang\underline{\BJ}\rang_{\GU\setminus\GO})]
         -(1-p)\Tr[\lang \underline{\BJ}^T \Tbb^\prime \underline{\BJ} \rang_{\GU\setminus\GO}]\right\},
\eeq{LB10a}
where $p$ is the volume fraction the body $\GO$ occupies within the cube $\GU$ (and the angular brackets $\lang\cdot\rang_{\GU\setminus\GO}$ 
denote volume averages over $\GU\setminus\GO$). This bound is not necessarily useful unless $\underline{\BJ}$ is close to being of the form \eq{LB1e}.
There is a wide variety of fluxes $\Bq$ for which $g$ is exactly computable. Let smooth potentials $\Ga$ and $\Gb$ be chosen in the neighborhood of the boundary of $\GO$ and set
\beq q_l=-n_k\left(J^0_{kl}+\frac{\Md^2\Ga}{\Md x_k \Md x_l}-\Gd_{kl}\Delta\Ga+\Ge_{klm}\frac{\Md\Gb}{\Md x_m}\right). 
\eeq{LB10b}
Then we can extend $\Ga$ and $\Gb$ to all of $\GU$ with periodic boundary conditions on $\p\GU$, and thus obtain a field $\underline{\BJ}$ given by \eq{LB1e} satisfying
 $\Bn \underline{\BJ} = -\Bq$ on $\p\Om$. This field attains equality in \eq{LB1d}:
\beq p\Tr[\lang \underline{\BJ}^T \Tbb^\prime \underline{\BJ} \rang]+(1-p)\Tr[\lang \underline{\BJ}^T \Tbb^\prime \underline{\BJ} \rang_{\GU\setminus\GO}=
\Tr[(\BJ^0)^T\Tbb^\prime\BJ^0].
\eeq{LB10bb}
The first term in this equation cannot be reduced, while maintaining the inequality $\eq{LB1d}$, by varying $\underline{\BJ}$ inside $\GO$ while keeping $\Bn \underline{\BJ} = -\Bq$ on $\p\Om$
(which ensures $\lang \underline{\BJ} \rang$ remains equal to $\BI$). Using \eq{LB1f} and integration by parts this establishes that
\beq g=\Tr[\lang \underline{\BJ}^T \Tbb^\prime \underline{\BJ} \rang]=\frac{1}{|\Om|} \int_{\p\Om} -[\Bx^T\Tbb^\prime\BJ^0+2\nabla\Ga]\cdot\Bq .
\eeq{LB10c}
In particular, if we take special Neumann conditions $\Bq=-\Bn$ on $\p\Om$, then \eq{LB10b} is satisfied with $\BJ^0=\BI$ and $\Ga=\Gb=0$, and the first equality in \eq{LB10c} implies $g=-3$.
It is an open question as to whether given any fluxes $q_1$, $q_2$ and $q_3$ (with no net flux through $\p\GO$) one can find $\Ga$, $\Gb$ and $\BJ^0$ such that \eq{LB10b} holds, however the 
following calculation suggests it might be possible. Suppose that there is a portion of $\p\GO$ which is flat with normal $\Bn=(0,0,1)^T$. Then on this flat section \eq{LB10b} implies
\beq q_1=-J^0_{31}-\frac{\Md\Ga_3}{\Md x_1}-\frac{\Md\Gb}{\Md x_2},\quad q_2=-J^0_{32}-\frac{\Md\Ga_3}{\Md x_2}+\frac{\Md\Gb}{\Md x_1},\quad q_3=-J^0_{33}+\frac{\Md^2\Ga}{\Md x_1^2}+\frac{\Md^2\Ga}{\Md x_2^2},
\eeq{LB10d}
where $\Ga_3=\Md\Ga_3/\Md x_3$. Thus possible values of $\Ga$, $\Ga_3$, and $\Gb$ on this flat interface may be found by solving Poisson's equations 
\beq \frac{\Md^2\Ga}{\Md x_1^2}+\frac{\Md^2\Ga}{\Md x_2^2}=q_3+J^0_{33}, \quad \frac{\Md^2\Ga_3}{\Md x_1^2}+\frac{\Md^2\Ga_3}{\Md x_2^2}=-\frac{\Md q_1}{\Md x_1}-\frac{\Md q_2}{\Md x_2}, \quad
\frac{\Md^2\Gb}{\Md x_1^2}+\frac{\Md^2\Gb}{\Md x_2^2}=\frac{\Md q_2}{\Md x_1}-\frac{\Md q_1}{\Md x_2}.
\eeq{LB10e}

Now we have from \eq{LB5}
 \begin{align}
  W_c^\prime & \ge \min_{\lang \underline{\BJ} \rang = \BI}
 \Tr \lang \underline{\BJ}^T [(\Gs^{-1}+c) \GL_h + (\Gs^{-1}-2c) \GL_s + \Gs^{-1}\GL_a ] \underline{\BJ} \rang \nonumber \\
 & \ge \min_{\lang \underline{\BJ} \rang = \BI}
 \Tr \lang \underline{\BJ}^T [(\Gs^{-1}+c) \GL_h ] \underline{\BJ} \rang, \label{LB11}
 \end{align}
by relaxing constraints. Observe that
 \begin{align*}
 \lang \underline{\BJ}^T [(\Gs^{-1}+c) \GL_h ] \underline{\BJ} \rang
 = \frac{1}{3} \lang (\Gs^{-1}+c) (\Tr \underline{\BJ})^2 \rang
 = \frac{\Gs_1^{-1}+c}{3|\Om|} \int_{\textrm{phase1}} (\Tr \underline{\BJ})^2 + \frac{\Gs_2^{-1}+c}{3|\Om|} \int_{\textrm{phase2}} (\Tr \underline{\BJ})^2.
 \end{align*}
It thus follows from Jensen's inequality that
 \beq
 \lang \underline{\BJ}^T [(\Gs^{-1}+c)\GL_h ] \underline{\BJ} \rang
 \ge \frac{\Gs_1^{-1}+c}{3f_1} \left(\frac{1}{|\Om|} \int_{\textrm{phase1}} \Tr \underline{\BJ} \right)^2 + \frac{\Gs_2^{-1}+c}{3f_2} \left(\frac{1}{|\Om|} \int_{\textrm{phase2}} \Tr \underline{\BJ} \right)^2 .
 \eeq{LB12}
We also have from the constraint $\lang \underline{\BJ} \rang = \BI$ that
 \beq
 \frac{1}{|\Om|} \int_{\textrm{phase1}} \Tr \underline{\BJ} + \frac{1}{|\Om|} \int_{\textrm{phase2}} \Tr \underline{\BJ} = 3.
 \eeq{LB13}
By minimizing the righthand side of \eq{LB12} under the constraint \eq{LB13} we get
 \beq
 W_c^\prime \ge \min_{\lang \underline{\BJ} \rang = \BI}
 \Tr \lang \underline{\BJ}^T [(\Gs^{-1}+c)\GL_h ] \underline{\BJ} \rang
 = \frac{3}{\frac{f_1}{(\Gs_1^{-1}+c)} + \frac{f_2}{(\Gs_2^{-1}+c)}},
 \eeq{LB14}
for all $c < \Gs_1^{-1}/2$. Letting $c \to \Gs_1^{-1}/2$, we have from \eq{LB10-4} that
 \beq
 \Tr \BA^\prime - \Gs_1^{-1} g^-/2 \ge \frac{3}{\frac{2f_1}{3\Gs_1^{-1}} + \frac{2f_2}{2\Gs_2^{-1}+\Gs_1}},
 \eeq{LB15}
which yields
 \beq
 f_1 =1-f_2\ge 1-\frac{2\Gs_1+\Gs_2}{2(\Gs_1 - \Gs_2)} \left[1-\frac{9}{2\Gs_1\Tr \BA^\prime -g^-}\right].
 \eeq{LB16}

With the special Neumann conditions $\Bq=-\Bn$ on $\p\Om$, we have $\lang \BJ \rang =\BI$ and
 \beq
 \BA^\prime = \lang \BJ^T \Gs^{-1} \BJ \rang = \lang \BJ^T \rang \lang \Gs^{-1} \BJ \rang = \BI \BGs_N^{-1} \lang \BJ \rang =\BGs_N^{-1},
 \eeq{LB17}
where $\BGs_N$ is the Neumann tensor. Since $g=-3$ in this case, we can take $g^-=-3$ giving
 \beq
 f_1 \ge 1-\frac{2\Gs_1+\Gs_2}{2(\Gs_1 - \Gs_2)} \left[1-\frac{9}{2\Gs_1\Tr(\BGs_N^{-1})+3}\right].
 \eeq{LB18}
As for the upper bound, a lower bound was obtained in \cite{Milton:2011:UBE}:
 \beq
 f_1 \ge 1-\frac{2\Gs_1+\Gs_2}{\Gs_1-\Gs_2} \frac{1}{\Gs_1 \Tr [(\Gs_1 \BI-\BGs_N)^{-1}] -1}.
 \eeq{LB19}
The bound in \eq{LB18} coincides with that in \eq{LB19} when $\BGs_N=\Gs_N\BI$. We could have recovered the bound \eq{LB19}, if we had defined
\beq 
W_c^\prime := \min_{\nabla \cdot \underline{\BJ} = 0, \ \lang \underline{\BJ} \rang = \BI, \atop \Bn \underline{\BJ} = -\Bq \ \textrm{on } \p\Om}
 \Tr \lang \BK^T\underline{\BJ}^T \Lbb_c^\prime(\underline{\BJ}\BK) \rang,\quad
 g:= \min_{\nabla \cdot \underline{\BJ} = 0, \ \lang \underline{\BJ} \rang = \BI, \atop \Bn \underline{\BJ} = -\Bq \ \textrm{on } \p\Om}
 \Tr \lang \BK^T \underline{\BJ}^T \Tbb^\prime(\underline{\BJ}\BK) \rang,
 \eeq{LB20}
and then optimized over the choice of $\BK$, in a similar manner as was done for the upper bound. This optimization is straightforward with
the Neumann conditions $\Bq=-\Bn$ on $\p\Om$, but in general might be difficult as it is not clear if the dependence of $g$ on $\BK$ is quadratic.
A way around this difficulty is to fix (independent of $\BK$) a $\underline{\BJ}$ in $\GU\setminus\GO$
such that $\underline{\BJ}$ satisfies periodic boundary conditions on $\p\GU$ and $\Bn \underline{\BJ} = -\Bq$ on $\p\Om$ and to take
\beq g^-  =  \frac{1}{p}\left\{\Tr[(p\BK+(1-p)\lang \underline{\BJ}\rang_{\GU\setminus\GO}\BK)^T\Tbb^\prime(p\BK+(1-p)\lang \underline{\BJ}\rang_{\GU\setminus\GO}\BK)] 
         -(1-p)\Tr[\lang (\underline{\BJ}\BK)^T \Tbb^\prime(\underline{\BJ}\BK) \rang_{\GU\setminus\GO}]\right\},
\eeq{LB21}
which then depends quadratically on $\BK$. However in the end one is still faced with the problem of choosing $\underline{\BJ}$ in $\GU\setminus\GO$ to
optimize the bound. The treatment given here can be generalized in other ways too, notably by allowing for translations $\Tbb^\prime$ which are not isotropic.

\section*{Acknowledgements}

HK is grateful for support from Ministry of Education, Sciences and Technology of Korea through NRF grants No. 2009-0085987 and 2010-0017532. 
GWM is grateful for support from the Mathematical Sciences Research Institute and from
National Science Foundation through grants DMS-070978 and DMS-***.

\bibliographystyle{siam}
\bibliography{/home/milton/tcbook,/home/milton/newref}

\end{document}